\documentclass[11pt]{article}
\usepackage{jahnel}

\title{Vanishing bond-percolation threshold implies\\ vanishing contact-process threshold}
\author{Benedikt Jahnel\,\orcidlink{0000-0002-4212-0065}\textsuperscript{1,2}\\
  \small\texttt{benedikt.jahnel@tu-braunschweig.de}
  \and Lukas L\"uchtrath\,\orcidlink{0000-0003-4969-806X}\textsuperscript{2}\\
  \small\texttt{lukas.luechtrath@wias-berlin.de}
  \and Christian M\"onch\,\orcidlink{0000-0002-6531-6482}\textsuperscript{3}\\
  \small\texttt{cmoench25@gmail.com}}
\date{\small
  \textsuperscript{1}Technische Universit\"at Braunschweig, Braunschweig, Germany\\
  \textsuperscript{2}Weierstrass Institute for Applied Analysis and Stochastics, Berlin, Germany\\
  \textsuperscript{3}Independent Researcher, Darmstadt, Germany\\[3mm]
  \normalsize\today}

\begin{document}
\maketitle

\begin{abstract}
We show that for any infinite, locally finite and connected graph $G$ and $p\in(0,1/2]$ the survival probability of the contact process with infection rate $\lambda\ge20p$ is at least $3/4$ times the bond percolation probability with retention probability $p$. In particular, $\lc(G)\le20\,\pc(G)$ whenever $\pc(G)<1/2$, and $\pc(G)=0$ implies $\lc(G)=0$. The proof uses a submartingale built from percolation clusters on finite balls. We contrast this with the SIR epidemic, in which every vertex is infected at most once, and for which there are graphs with $\pc(G)=0$ on which the SIR epidemic dies out almost surely for every infection rate.

\medskip\noindent\emph{This note is not intended for journal publication.}
\end{abstract}

\noindent\textbf{AMS-MSC 2020}: 60K35 (primary); 82B43 (secondary)


\section{Introduction}\label{sec:intro}

The contact process is a basic model for epidemic spread through a structured population and one of the simplest interacting particle systems that exhibits a survival/extinction phase transition~\cite{Harris1974,Liggett1985,Liggett1999}. Throughout, $G$ is an infinite, locally finite and connected graph with root $o$. In the contact process on $G$, infected vertices infect each of their neighbours at rate $\lambda>0$ and recover at rate $1$. Starting from the single infected vertex $o$, the process survives with positive probability if $\lambda$ exceeds the critical infection rate $\lc(G)$, and it dies out almost surely if $\lambda<\lc(G)$. Since $G$ contains a ray and $\lc$ is monotone under taking subgraphs, we have $\lc(G)\in[0,\lc(\mathbb N)]$, where $\lc(\mathbb N)<\infty$ is the critical value on the half-line~\cite{Harris1974,Liggett1999}. Note that on graphs that expand sufficiently strongly, the contact process might infect infinitely many vertices in finite time, which we treat as part of the survival phase.

Recall that Bernoulli bond percolation with parameter $p\in[0,1]$ on $G$ has an infinite cluster with positive probability if $p$ exceeds the critical value $\pc(G)\in[0,1]$, and almost surely no infinite cluster if $p<\pc(G)$.

Both critical values are positive on graphs with uniformly bounded degrees, by comparison with branching processes. On graphs with unbounded degrees, both critical values may vanish. In percolation theory, this robustness phenomenon occurs, for example, on Galton--Watson trees whose offspring distribution has infinite mean~\cite{Lyons1990}, in scale-free percolation with infinite degree variance~\cite{DeijfenHofstadHooghiemstra2013}, and in weight-dependent random connection models with sufficiently strong preference for high-weight vertices~\cite{GracarLuechtrathMoerters2021}. For the contact process, survival for every $\lambda>0$ is known on Galton--Watson trees with sufficiently heavy-tailed offspring distributions, conditional on the tree being infinite~\cite{Pemantle1992,HuangDurrett2020,BarthaKomjathyValesin2026}. A closely related phenomenon is the persistence for stretched-exponential time of the infection on locally tree-like random graphs~\cite{ChatterjeeDurrett2009}. The underlying mechanism is well understood: a vertex of degree $d$ keeps the infection in its neighbourhood for a time of order $\exp(c\lambda^2d)$, which is enough to transmit the infection to the next vertex of high degree. In our work~\cite[Section~1.1]{JahnelLuechtrathMoench2025} we mentioned that $\pc(G)=0$ should imply $\lc(G)=0$.

This note is motivated by a question of Daniel Valesin regarding this claim in an earlier version of~\cite{JahnelLuechtrathMoench2025}, and its purpose is to provide a rigorous proof of that implication. We write $\theta_{\mathsf{CP}}(\lambda)$ for the probability that the contact process with infection rate $\lambda$ started from $\{o\}$ survives, and $\theta_{\mathsf{bond}}(p)$ for the probability that the open cluster of $o$ in Bernoulli bond percolation with parameter $p$ is infinite. Our main result is the following comparison.

\begin{theorem}\label{thm:comparison}
Let $G$ be an infinite, locally finite and connected graph with root $o$. If $0<p\le1/2$ and $\lambda\ge20p$, then
\[
\theta_{\mathsf{CP}}(\lambda)\ge\tfrac34\,\theta_{\mathsf{bond}}(p).
\]
\end{theorem}

Since $\theta_{\mathsf{bond}}(p)>0$ for every $p>\pc(G)$, \cref{thm:comparison} immediately yields a bound on the critical infection rate in terms of the critical percolation parameter.

\begin{corollary}\label{cor:comparison}
Let $G$ be an infinite, locally finite and connected graph with $\pc(G)<1/2$. Then $\lc(G)\le20\,\pc(G)$. In particular, $\pc(G)=0$ implies $\lc(G)=0$.
\end{corollary}

\begin{proof}
For every $p\in(\pc(G),1/2]$ we have $\theta_{\mathsf{bond}}(p)>0$, and hence $\theta_{\mathsf{CP}}(20p)>0$ by \cref{thm:comparison}. Thus $\lc(G)\le20p$, and letting $p\downarrow\pc(G)$ proves the claim.
\end{proof}

The proof of \cref{thm:comparison} in \cref{sec:proof-comparison} has three steps. Given a finite graph $H$ and a non-empty set $B$ of its vertices that the infection should reach, we consider the function $F(A)=1-\mathbb E_p[4^{-\#(A\cap\mathcal C_B)}]$ on sets of infected vertices, where $\mathcal C_B$ is the set of vertices connected to $B$ by open paths. A weighted estimate on the boundary of percolation clusters shows that $F$ is subharmonic for the contact process stopped upon hitting $B$, as long as $\lambda\ge20p$. Applied to balls and spheres, this shows that the contact process reaches distance $n$ with probability at least $(3/4)\,\theta_{\mathsf{bond}}(p)$. A short exhaustion argument then turns spatial spread into survival in time.

\paragraph{Outline.} In \cref{sec:results} we introduce the setting and state further results: a bound on trees that is sharp to first order on regular trees, an example showing that a comparison with the SIR epidemic fails on general graphs, and a remark on the converse of \cref{cor:comparison}. \Cref{sec:proof-comparison} contains the proof of \cref{thm:comparison}, \cref{sec:proof-tree} the proof of the bound on trees and \cref{sec:proof-chain} the proof concerning the SIR example.

\section{Setting and further results}\label{sec:results}

\subsection{Contact process and percolation}\label{sec:setting}

Throughout, $G=(V,E)$ is an infinite, locally finite and connected simple graph with a distinguished root $o\in V$. We write $x\sim y$ if $\{x,y\}\in E$, $N(x)\coloneqq\{y\in V\colon y\sim x\}$ and $\deg(x)\coloneqq\#N(x)$. The graph distance is denoted by $\dist$.

\begin{figure}[ht]
\centering
\begin{tikzpicture}[xscale=1.25,yscale=0.85,
  rec/.style={font=\small},
  arr/.style={-{Stealth[length=2mm]},black!60,line width=0.6pt},
  used/.style={-{Stealth[length=2.2mm]},red!75!black,line width=1.2pt},
  inf/.style={red!75!black,line width=2.2pt}]
\foreach \c/\lab in {0/{-2},1/{-1},2/{0},3/{1},4/{2}} {
  \draw[black!30,line width=0.6pt] (\c,0) -- (\c,5);
  \node[below] at (\c,0) {$\lab$};
}
\draw[-{Stealth},black!50] (-0.6,0) -- (-0.6,5.3) node[above] {time};
\draw[inf] (2,0) -- (2,2.0);
\draw[inf] (3,0.8) -- (3,1.2);
\draw[inf] (4,1.0) -- (4,3.0);
\draw[inf] (1,1.5) -- (1,3.2);
\draw[inf] (0,2.2) -- (0,2.6);
\draw[inf] (2,2.4) -- (2,5);
\draw[inf] (3,3.5) -- (3,4.3);
\draw[arr] (4,0.5) -- (3,0.5);
\draw[arr] (3,2.7) -- (2,2.7);
\draw[arr] (0,4.0) -- (1,4.0);
\draw[used] (2,0.8) -- (3,0.8);
\draw[used] (3,1.0) -- (4,1.0);
\draw[used] (2,1.5) -- (1,1.5);
\draw[used] (1,2.2) -- (0,2.2);
\draw[used] (1,2.4) -- (2,2.4);
\draw[used] (2,3.5) -- (3,3.5);
\foreach \c/\t in {2/2.0,3/1.2,3/4.3,1/3.2,0/2.6,4/3.0,0/1.0,4/4.6} \node[rec] at (\c,\t) {$\times$};
\node[right] at (4.2,5) {$t$};
\draw[black!40,dashed] (-0.3,5) -- (4.3,5);
\end{tikzpicture}
\caption{Harris' graphical representation on a part of $\mathbb Z$ with $o=0$. Recovery marks are drawn as crosses and infection arrows as horizontal arrows. Infected space-time points are drawn in red; arrows emanating from healthy vertices (grey) have no effect. The vertex $o$ recovers, is reinfected by its left neighbour, and is the only infected vertex at time $t$.}
\label{fig:harris}
\end{figure}
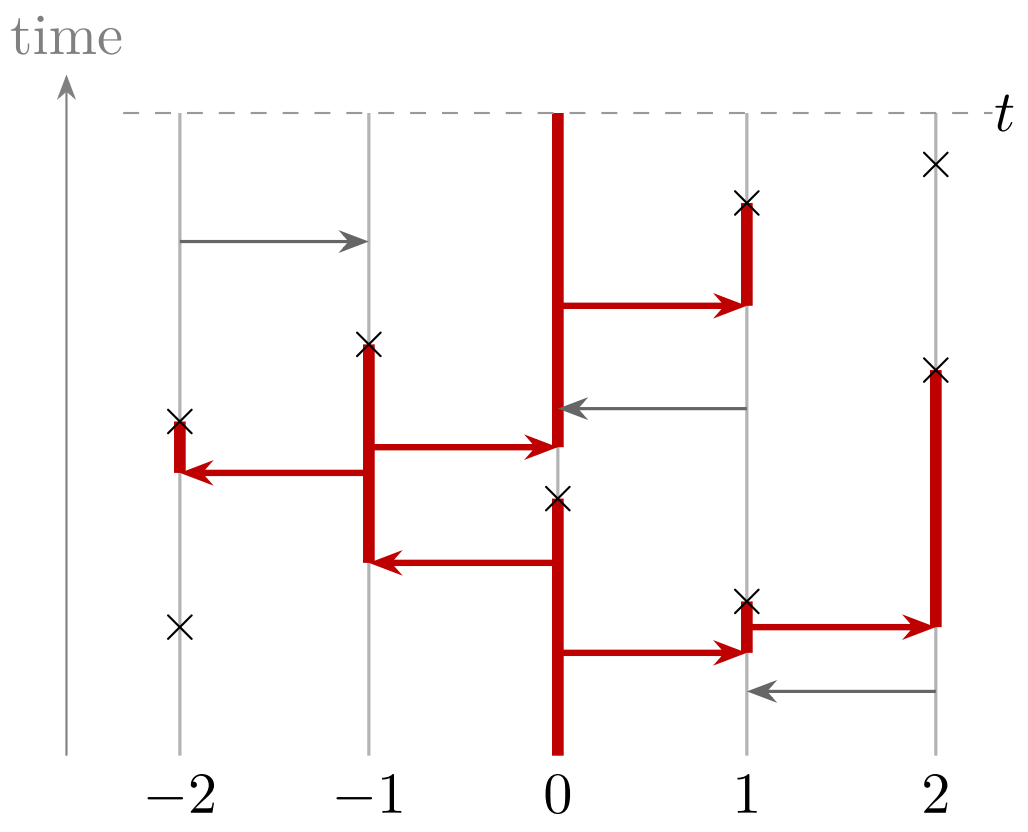

We construct the contact process from Harris' graphical representation~\cite{Harris1974,Liggett1985}. For $\lambda>0$ let $D_x$, $x\in V$, be Poisson point processes on $[0,\infty)$ with intensity $1$, the \emph{recovery marks}, and let $A_{(x,y)}$, $x\sim y$, be Poisson point processes with intensity $\lambda$, the \emph{infection arrows} from $x$ to $y$, all independent under $\mathbb P_\lambda$. An \emph{infection path} from $(x,s)$ to $(y,t)$ is a sequence $x=x_0,\dots,x_j=y$ with times $s=t_0<t_1<\dots<t_j\le t_{j+1}=t$ such that $t_i\in A_{(x_{i-1},x_i)}$ for $1\le i\le j$, $D_{x_0}\cap[t_0,t_1]=\emptyset$ and $D_{x_i}\cap(t_i,t_{i+1}]=\emptyset$ for $1\le i\le j$; see \cref{fig:harris}. For $A\subset V$ let $\xi_{s,t}(A)$ be the set of $y\in V$ such that there is an infection path from $(x,s)$ to $(y,t)$ for some $x\in A$, and write $\xi_t(A)\coloneqq\xi_{0,t}(A)$ and $\xi_t\coloneqq\xi_t(\{o\})$. As long as $\xi_t$ is finite, this is the contact process with the generator given in~\cite[Section~1.1]{JahnelLuechtrathMoench2025}. If the degrees grow very fast, $\xi_t$ may become infinite in finite time; the graphical representation defines $\xi_t$ for all $t\ge0$ also in this case, and all statements refer to this definition. On a subgraph $H$ of $G$ containing $o$ we use the restriction of the graphical representation, and we write $\xi^{H}_t$ for the resulting contact process on $H$ started from $\{o\}$.

We write $\cF_t$ for the $\sigma$-algebra generated by the restrictions of all processes $D_x$ and $A_{(x,y)}$ to $[0,t]$. By construction, $\xi_{s,t+u}(A)=\xi_{t,t+u}(\xi_{s,t}(A))$ for $s\le t$ and $u\ge0$, and the restrictions of the Poisson processes to $(t,\infty)$ are independent of $\cF_t$. We use this Markov property at deterministic times and at stopping times with countably many values. On finite graphs, we also use the strong Markov property at arbitrary stopping times. The process \emph{survives} if $\xi_t\neq\emptyset$ for all $t\ge0$, and we define
\[
\theta_{\mathsf{CP}}(\lambda)\coloneqq\mathbb P_\lambda(\xi_t\neq\emptyset\text{ for all }t\ge0),\qquad\lc(G)\coloneqq\inf\{\lambda>0\colon\theta_{\mathsf{CP}}(\lambda)>0\}.
\]
It is immediate from the graphical representation that the contact process is monotone in the underlying graph and in the initial configuration, both in the sense of set inclusion: if $A\subset B\subset V$, then $\xi_t(A)\subset\xi_t(B)$, and if $H$ is a subgraph of $G$, then the contact process on $H$ constructed from the restriction of the graphical representation is contained in the one on $G$. Monotonicity in the initial configuration is also referred to as attractiveness. Similarly, the process is monotone in $\lambda$ by the standard coupling of the graphical representations. In particular, $\lc(G)$ is the unique threshold between extinction and survival, it does not depend on the choice of the root since $G$ is connected, and $\lc(G)\le\lc(H)$ for every infinite connected subgraph $H$ of $G$.

Let $\mathbb P_p$ denote Bernoulli bond percolation on $G$ with parameter $p\in[0,1]$, let $\theta_{\mathsf{bond}}(p)\coloneqq\mathbb P_p(o\leftrightarrow\infty)$, where $\{o\leftrightarrow\infty\}$ is the event that the open cluster of $o$ is infinite, and let $\pc(G)\coloneqq\sup\{p\in[0,1]\colon\theta_{\mathsf{bond}}(p)=0\}$. Also $\pc(G)$ does not depend on the choice of the root, and $\pc(G)\le\pc(H)$ for every infinite connected subgraph $H$ of $G$.

\subsection{Trees}\label{sec:trees}

Let $T$ be an infinite, locally finite tree with root $o$. For a vertex $x$ we write $|x|\coloneqq\dist(o,x)$, and for an edge $e$ we write $|e|$ for the larger of the distances of its endpoints to $o$. A \emph{cutset} is a finite set $\Pi$ of edges such that every infinite self-avoiding path starting at $o$ uses an edge of $\Pi$. The \emph{branching number} of $T$ is
\[
\br(T)\coloneqq\sup\Bigl\{\rho\ge1\colon\inf_{\Pi}\sum_{e\in\Pi}\rho^{-|e|}>0\Bigr\},
\]
where the infimum runs over all cutsets; it does not depend on the choice of the root~\cite[Chapter~1]{LyonsPeres2016}. Lyons~\cite{Lyons1990} proved that
\begin{equation}\label{eq:lyons}
\pc(T)=1/\br(T).
\end{equation}

\begin{proposition}
\label{prop:tree}
Let $T$ be an infinite, locally finite tree with $\pc(T)<1$. Then
\[
\lc(T)\le\frac{\pc(T)}{1-\pc(T)}=\frac1{\br(T)-1}.
\]
In particular, $\lc(T)\le2\,\pc(T)$ if $\pc(T)\le1/2$.
\end{proposition}

The identity in \cref{prop:tree} is~\eqref{eq:lyons}. Recall from \cref{sec:intro} that also $\lc(T)\le\lc(\mathbb N)$, since $T$ contains a ray; this bound is better if $\pc(T)$ is close to $1$. The proof of \cref{prop:tree} in \cref{sec:proof-tree} compares the contact process with the SIR epidemic of \cref{sec:chain}, which on trees is a quasi-independent percolation in the sense of Lyons~\cite{Lyons1989}, and uses the exhaustion argument of \cref{sec:exhaustion} to turn spatial spread into survival.

\begin{remark}\label{rem:regular-tree}
The bound of \cref{prop:tree} is sharp to first order on regular trees. Let $T_d$ be the tree in which every vertex has $d\ge2$ children, so that $\pc(T_d)=1/d$ and \cref{prop:tree} gives $\lc(T_d)\le1/(d-1)$. Conversely, on every graph with maximal degree $\Delta$, every infected vertex creates new infections at rate at most $\lambda\Delta$ and recovers at rate $1$. Comparison with a branching random walk yields $\mathbb E_\lambda[\#\xi_t]\le\e^{(\lambda\Delta-1)t}$, and hence extinction if $\lambda<1/\Delta$. Since $T_d$ has maximal degree $d+1$, we obtain the classical bounds, see~\cite[Part~I]{Liggett1999},
\[
\frac1{d+1}\le\lc(T_d)\le\frac1{d-1}.
\]
\end{remark}

By monotonicity, the bound of \cref{prop:tree} transfers to every graph $G$ containing $T$ as a subgraph. For a general graph $G$ let
\[
\pc^{\mathsf{tree}}(G)\coloneqq\inf\{\pc(T)\colon T\text{ is an infinite subtree of }G\}\ \ge\ \pc(G).
\]
Here and below, we interpret $u/(1-u)$ as $+\infty$ when $u=1$. Then $\lc(G)\le\pc^{\mathsf{tree}}(G)/(1-\pc^{\mathsf{tree}}(G))\wedge\lc(\mathbb N)$, since $u\mapsto u/(1-u)$ is continuous and increasing. This improves the constant of \cref{cor:comparison} for graphs with $\pc^{\mathsf{tree}}(G)$ close to $\pc(G)$, but \cref{prop:chain} below shows that $\pc^{\mathsf{tree}}(G)$ can be much larger than $\pc(G)$.

\subsection{The converse implication}\label{sec:converse}

\begin{figure}[ht]
\centering
\begin{tikzpicture}[scale=1,
  hub/.style={circle,fill=black,inner sep=2pt},
  leaf/.style={circle,draw=black,fill=white,inner sep=1.1pt}]
\draw (0,0) -- (8.6,0);
\node at (9.1,0) {$\cdots$};
\foreach \k/\x in {1/0.3,2/1.9,3/4.3} {
  \node[hub,label=below:{$\k$}] (h\k) at (\x,0) {};
  \pgfmathtruncatemacro{\n}{\k*\k}
  \foreach \i in {1,...,\n} {
    \pgfmathsetmacro{\ang}{180*\i/(\n+1)}
    \pgfmathsetmacro{\r}{0.55+0.22*\k}
    \draw (h\k) -- ++(\ang:\r) node[leaf] {};
  }
}
\node[hub,label=below:{$4$}] (h4) at (7.6,0) {};
\foreach \i in {1,...,16} {
  \pgfmathsetmacro{\ang}{180*\i/17}
  \draw (h4) -- ++(\ang:1.45) node[leaf] {};
}
\end{tikzpicture}
\caption{The graph of \cref{sec:converse}: the half-line with $k^2$ leaves attached to the vertex $k$. Every edge of the half-line is a cut edge, so that $\pc(G)=1$, while the stars of growing size yield $\lc(G)=0$.}
\label{fig:stars}
\end{figure}
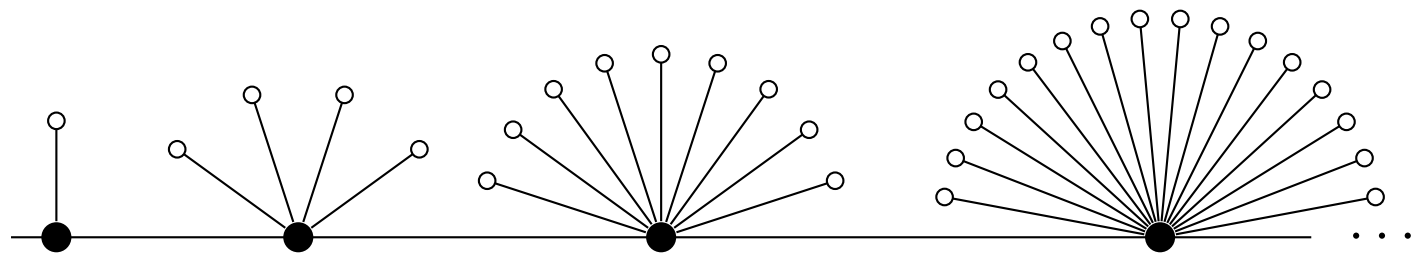

The converse of \cref{cor:comparison} fails. Bartha, Komj\'athy and Valesin~\cite[Definition~6.3]{BarthaKomjathyValesin2026} consider the graph $H_{K,\ell}$, a one-ended row of stars with $K$ leaves each, whose centres are joined by paths with $\ell$ inner vertices. By~\cite[Lemma~6.8]{BarthaKomjathyValesin2026}, for every $\lambda>0$ and every $\ell\colon\mathbb N\to\mathbb N$ with $\ell(K)=o(K)$, the contact process on $H_{K,\ell(K)}$ survives, even locally, with positive probability for all sufficiently large $K$. We use the $\mu=0$ case of their result for degree-penalised contact processes with infection rate $\lambda/(\deg(x)\deg(y))^\mu$; for general $\mu\in[0,1/2)$, their hypothesis is $\ell(K)=o(K^{1-2\mu})$. The proof uses that a star with $K$ leaves keeps the infection alive for a time $\exp(c\lambda^2K)$~\cite[Claim~6.6]{BarthaKomjathyValesin2026}, see also~\cite[Lemma~3.1]{MountfordValesinYao2013}, that it passes the infection along a path of length $o(K)$ to the next star~\cite[Claim~6.7]{BarthaKomjathyValesin2026}, and a renormalisation argument. By monotonicity, $\lc(G)=0$ for every graph $G$ containing a copy of $H_{K,\ell(K)}$ for all sufficiently large $K$. For example, let $G$ be the half-line $\mathbb N$ with $k^2$ additional leaves attached to every $k\in\mathbb N$, see \cref{fig:stars}. For every $K$, the vertices $2j$ with $4j^2\ge K$ as centres and the vertices $2j+1$ as connecting paths yield a copy of $H_{K,1}$, so that $\lc(G)=0$. On the other hand, every edge $\{k,k+1\}$ is a cut edge, and hence $\pc(G)=1$.

\subsection{The single-infection process}\label{sec:chain}

Let the random variables $T_x$, $x\in V$, be independent and exponentially distributed with parameter~$1$. Let the random variables $E_{(x,y)}$, $x\sim y$, be independent of them, independent of each other and exponentially distributed with parameter~$\lambda$. The directed edge $(x,y)$ is \emph{open} if $E_{(x,y)}<T_x$, and we let $\cS_\lambda$ be the set of vertices that can be reached from $o$ by a directed path of open edges. The set $\cS_\lambda$ is the set of vertices that are ever infected in the SIR epidemic started from $o$ in which infected vertices infect their neighbours at rate $\lambda$ and are removed after an exponential time with parameter $1$; in this epidemic every vertex is infected at most once~\cite{Kuulasmaa1982}. Every directed edge is open with probability $\lambda/(1+\lambda)$, but edges leaving the same vertex are positively correlated.

\begin{figure}[ht]
\centering
\begin{tikzpicture}[scale=1,
  hub/.style={circle,fill=black,inner sep=2pt},
  mid/.style={circle,draw=black,fill=white,inner sep=1.3pt}]
\foreach \k/\x in {1/0,2/2.4,3/4.8,4/7.2} \node[hub,label=below:{$v_{\k}$}] (v\k) at (\x,0) {};
\foreach \k/\x in {1/0,2/2.4,3/4.8} {
  \pgfmathtruncatemacro{\kk}{\k+1}
  \foreach \i in {1,...,\k} {
    \pgfmathsetmacro{\y}{(\i-(\k+1)/2)*0.55}
    \draw (v\k) -- (\x+1.2,\y) -- (v\kk);
    \node[mid] at (\x+1.2,\y) {};
  }
}
\foreach \i in {1,...,4} {
  \pgfmathsetmacro{\y}{(\i-2.5)*0.55}
  \draw (v4) -- (8.4,\y);
  \node[mid] at (8.4,\y) {};
}
\node at (9.2,0) {$\cdots$};
\end{tikzpicture}
\caption{The graph $G_{\mathsf{chain}}$ of \cref{prop:chain}: the vertices $v_k$ and $v_{k+1}$ are connected by $k$ disjoint paths of length $2$ through the vertices $w_{k,1},\dots,w_{k,k}$ (white).}
\label{fig:chain}
\end{figure}
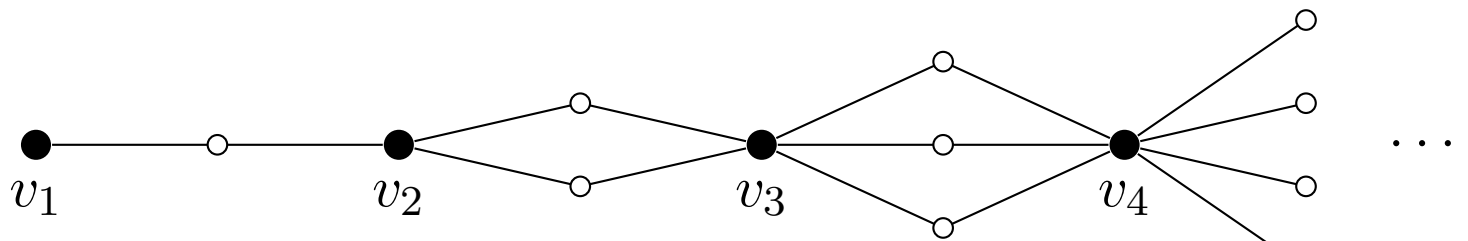

Define the graph $G_{\mathsf{chain}}$ as follows. Its vertices are $v_k$ for $k\in\mathbb N$ and $w_{k,i}$ for $k\in\mathbb N$ and $1\le i\le k$, and its edges are $\{v_k,w_{k,i}\}$ and $\{w_{k,i},v_{k+1}\}$ for $k\in\mathbb N$ and $1\le i\le k$. Thus $v_k$ and $v_{k+1}$ are connected by $k$ disjoint paths of length $2$, see \cref{fig:chain}. The root is $o\coloneqq v_1$.

\newpage
\begin{proposition}\label{prop:chain}
The graph $G_{\mathsf{chain}}$ satisfies
\begin{enumerate}[label=(\alph*)]
\item\label{it:chain-perc} $\pc(G_{\mathsf{chain}})=0$, while $\pc^{\mathsf{tree}}(G_{\mathsf{chain}})=1$;
\item\label{it:chain-sir} $\#\cS_\lambda<\infty$ almost surely for every $\lambda>0$;
\item\label{it:chain-cp} $\lc(G_{\mathsf{chain}})=0$.
\end{enumerate}
\end{proposition}

Part~\ref{it:chain-sir} shows that no comparison between Bernoulli bond percolation and the single-infection process can yield $\lc(G)=0$ from $\pc(G)=0$ in general. Part~\ref{it:chain-cp} follows from \cref{cor:comparison} and part~\ref{it:chain-perc}, and also from \cref{sec:converse}; survival is due to reinfections of the vertices $v_k$, which are absent in the single-infection process. The same phenomenon occurs if $v_k$ and $v_{k+1}$ are connected by only $\lceil\log^2(k+1)\rceil$ disjoint paths of length $2$: bond percolation still crosses all blocks with positive probability, while the events that $v_k$ is removed before it transmits along any of its forward edges are independent with probabilities of order $1/(\lambda\log^2k)$, whose sum diverges. On trees, by contrast, the single-infection process is a quasi-independent percolation with the same edge marginals as Bernoulli bond percolation, and this is the basis of the proof of \cref{prop:tree}.

\begin{question}\label{qu:constant}
Does $\lc(G)\le(1+o(1))\,\pc(G)$ hold as $\pc(G)\to0$ for all infinite, locally finite and connected graphs $G$, or even $\lc(G)\le\pc(G)/(1-\pc(G))$ for every such graph? \Cref{cor:comparison} gives the first statement with the constant $20$ in place of $1+o(1)$, and \cref{prop:tree} proves the second one for trees. By \cref{rem:regular-tree}, the constant $1$ cannot be replaced by a smaller one.
\end{question}

\section{Proof of \cref*{thm:comparison}}\label{sec:proof-comparison}

Let $H$ be a finite simple graph with vertex set $V_H$ and let $\emptyset\neq B\subset V_H$. For a percolation configuration on $H$ let $\mathcal C_B$ be the set of vertices of $H$ that are connected to $B$ by an open path; in particular, $B\subset\mathcal C_B$. The set $\mathcal C_B$ does not depend on any set of infected vertices.

\subsection{A weighted boundary estimate}\label{sec:boundary}

\begin{figure}[ht]
\centering
\begin{tikzpicture}[scale=1.05,
  vtx/.style={circle,fill=black,inner sep=1.3pt},
  inf/.style={circle,fill=red!75!black,inner sep=2.1pt},
  bnd/.style={rectangle,fill=black,inner sep=2pt},
  open/.style={blue!70!black,line width=1.3pt},
  closed/.style={black!35,line width=0.6pt},
  mpair/.style={-{Stealth[length=2.2mm]},red!75!black,dashed,line width=0.8pt}]
\draw[black!25,dashed] (0,0) circle (2.8);
\coordinate (o) at (0,0);
\coordinate (a) at (1,0.6);   \coordinate (b) at (1.9,1.1);  \coordinate (s1) at (2.42,1.4);
\coordinate (c) at (-0.9,0.8);\coordinate (d) at (-1.7,1.5); \coordinate (s2) at (-1.98,1.98);
\coordinate (e) at (0.3,-1.1);\coordinate (f) at (0.9,-1.9); \coordinate (s3) at (1.4,-2.42);
\coordinate (g) at (-1.2,-0.7);\coordinate (h) at (-2.0,-1.2);\coordinate (s4) at (-2.42,-1.4);
\coordinate (s5) at (2.8,0); \coordinate (s6) at (0,2.8); \coordinate (s7) at (-2.8,0); \coordinate (s8) at (0,-2.8);
\draw[closed] (o)--(c) (c)--(d) (o)--(e) (o)--(g) (g)--(h) (h)--(s4) (b)--(s5) (d)--(s6) (h)--(s7) (f)--(s8) (a)--(e);
\draw[open] (o)--(a) (a)--(b) (b)--(s1) (e)--(f) (f)--(s3) (d)--(s2);
\draw[mpair] (a) to[bend left=25] (b);
\draw[mpair] (e) to[bend right=25] (f);
\draw[mpair] (c) to[bend left=25] (d);
\foreach \v in {b,d,f,h} \node[vtx] at (\v) {};
\foreach \v in {s1,s2,s3,s4,s5,s6,s7,s8} \node[bnd] at (\v) {};
\foreach \v in {o,a,c,e,g} \node[inf] at (\v) {};
\node[below left=1pt] at (o) {$o$};
\node[above left=0pt] at (a) {$x_1$};
\node[left=1pt] at (e) {$x_2$};
\node[below=2pt] at (c) {$x_3$};
\node[above=1pt] at (g) {$x_4$};
\node at (3.35,2.2) {$B$};
\end{tikzpicture}
\caption{Illustration of \cref{lem:boundary} for a finite graph $H$ and given $B$. Open edges are drawn in blue and closed edges in grey, and $\mathcal C_B$ is the set of vertices joined to $B$ by blue paths. The set $A=\{o,x_1,x_2,x_3,x_4\}$ is drawn in red. Here $A\cap\mathcal C_B=\{o,x_1,x_2\}$, so that $\mathcal N=3$, and the three dashed arrows mark the pairs counted by $\mathcal M$; the pair from $x_3$ is counted although the corresponding edge is closed.}
\label{fig:boundary}
\end{figure}
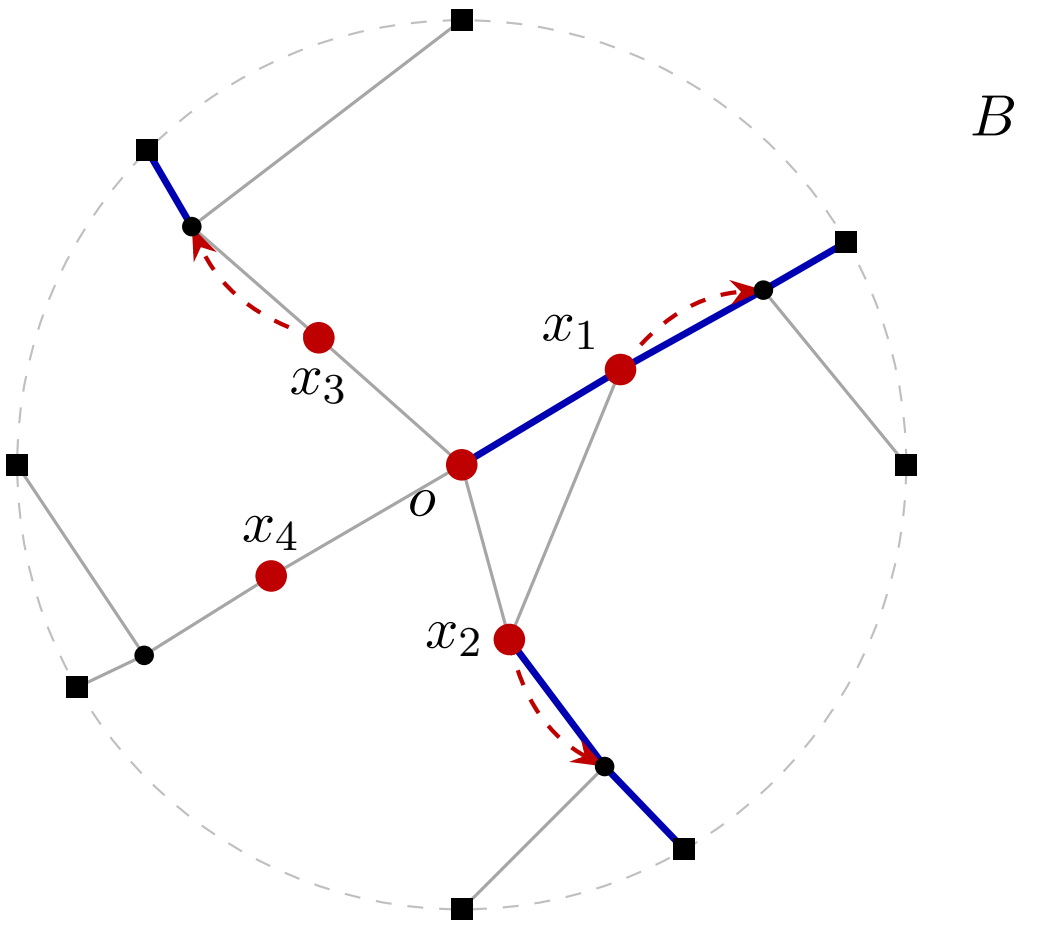

\begin{lemma}\label{lem:boundary}
Let $0<p\le1/2$ and $A\subset V_H\setminus B$. Define
\[
\mathcal N\coloneqq\#(A\cap\mathcal C_B),\qquad\mathcal M\coloneqq\#\bigl\{(x,y)\colon x\in A,\ y\in\mathcal C_B\setminus A,\ x\sim y\bigr\},
\]
where the edges counted by $\mathcal M$ need not be open; see \cref{fig:boundary}. Then
\begin{equation}\label{eq:boundary}
\mathbb E_p\bigl[\mathcal N\,4^{-\mathcal N}\bigr]\le5p\,\mathbb E_p\bigl[\mathcal M\,4^{-\mathcal N}\bigr].
\end{equation}
\end{lemma}

\begin{figure}[ht]
\centering
\begin{tikzpicture}[scale=1,
  vtx/.style={circle,fill=black,inner sep=1.2pt},
  open/.style={blue!70!black,line width=1.2pt},
  mm/.style={red!75!black,dashed,line width=0.9pt}]
\draw[black!60,rounded corners=12pt] (-3.2,-1.7) rectangle (0.2,1.7);
\node at (-2.85,1.4) {$A$};
\draw[fill=black!8,rounded corners=8pt] (-1.4,-0.7) rectangle (-0.1,0.7);
\node at (-0.75,0.95) {$K$};
\node[vtx] (k1) at (-1.05,0.3) {}; \node[vtx] (k2) at (-0.45,-0.25) {}; \node[vtx] (k3) at (-1.05,-0.4) {};
\draw[open] (k1)--(k2)--(k3)--(k1);
\draw[black!35,rounded corners=6pt] (-2.9,0.3) rectangle (-2.0,1.2);
\draw[black!35,rounded corners=6pt] (-2.9,-1.4) rectangle (-1.9,-0.3);
\draw[blue!70!black,line width=0.8pt,rounded corners=14pt] (1.4,-1.9) rectangle (5.4,1.9);
\node[blue!70!black] at (4.9,1.55) {$\mathcal C_0$};
\draw[line width=2pt] (5.4,-1.5) -- (5.4,1.5);
\node[right] at (5.5,0) {$B$};
\node[vtx] (y1) at (1.9,0.8) {}; \node[vtx] (y2) at (2.0,-0.2) {}; \node[vtx] (y3) at (1.9,-1.2) {};
\node[vtx] (z1) at (1.0,2.3) {}; \node[vtx] (z2) at (0.9,-2.3) {};
\draw[mm] (k1)--(y1); \draw[mm] (k2)--(y2); \draw[mm] (k3)--(y3);
\draw[black!35] (k1)--(z1); \draw[black!35] (k3)--(z2);
\draw[open] (y1)--(3.2,1.0)--(5.4,1.1);
\draw[open] (y2)--(3.4,-0.1)--(5.4,-0.3);
\draw[open] (y3)--(3.0,-1.3)--(4.2,-0.9)--(5.4,-1.0);
\end{tikzpicture}
\caption{The decomposition in the proof of \cref{lem:boundary}. The open component $K$ of $A$ is connected to $B$ if and only if one of the $m$ dashed edges from $K$ into $\mathcal C_0\setminus A$ is open, where $\mathcal C_0$ is defined without $K$. Given the revealed configuration, edges from $K$ to vertices outside $\mathcal C_0$ (grey) and edges to other components of $A$ are irrelevant to whether $K$ connects to $B$.}
\label{fig:decomposition}
\end{figure}
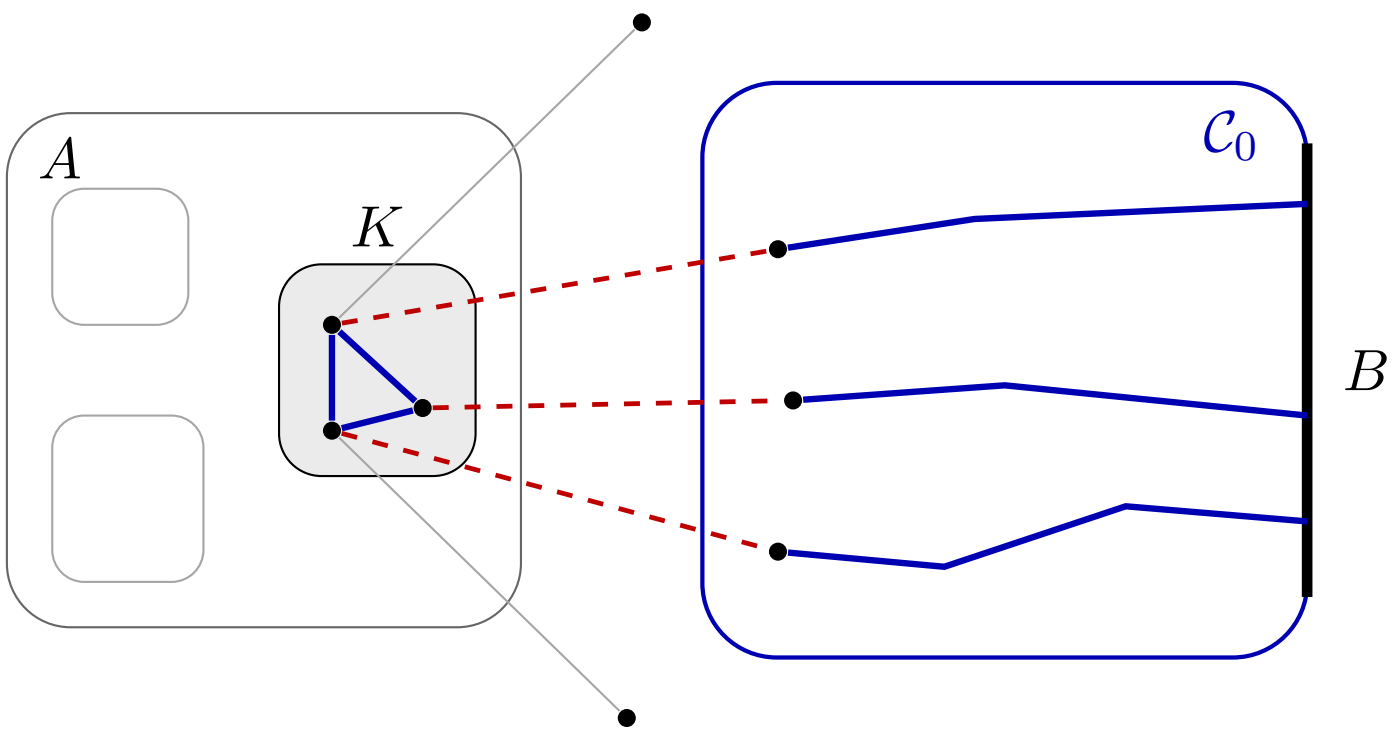

\begin{proof}
If $A=\emptyset$, both sides vanish. Otherwise, we condition on the states of all edges with both endpoints in $A$ and denote the open components of $A$ with respect to these edges by $K_1,\dots,K_r$, with $k_i\coloneqq\#K_i\ge1$. Every $K_i$ is connected by open edges, hence either $K_i\subset\mathcal C_B$ or $K_i\cap\mathcal C_B=\emptyset$. With
\[
\mathcal M_i\coloneqq\#\bigl\{(x,y)\colon x\in K_i,\ y\in\mathcal C_B\setminus A,\ x\sim y\bigr\}
\]
we therefore have $\mathcal M=\sum_i\mathcal M_i$ and $\mathcal N=\sum_ik_i\1\{K_i\subset\mathcal C_B\}$. It suffices to show that, for every $i$,
\begin{equation}\label{eq:boundary-i}
k_i\,\mathbb E_p\bigl[4^{-\mathcal N}\1\{K_i\subset\mathcal C_B\}\bigm\lvert\mathcal R_0\bigr]\le5p\,\mathbb E_p\bigl[4^{-\mathcal N}\mathcal M_i\bigm\lvert\mathcal R_0\bigr],
\end{equation}
where $\mathcal R_0$ denotes the $\sigma$-algebra generated by the states of the edges inside $A$. Summing~\eqref{eq:boundary-i} over $i$ and taking expectations yields~\eqref{eq:boundary}.

Fix $i$ and write $K\coloneqq K_i$ and $k\coloneqq k_i$. Let $\mathcal R\supset\mathcal R_0$ be the $\sigma$-algebra generated by the states of all edges except those between $K$ and $V_H\setminus A$. Let $\mathcal C_0$ be the set of vertices connected to $B$ by an open path in the graph obtained from $H$ by deleting $K$, see \cref{fig:decomposition}, and define
\[
N_0\coloneqq\#\bigl((A\setminus K)\cap\mathcal C_0\bigr),\qquad m\coloneqq\#\bigl\{(x,y)\colon x\in K,\ y\in\mathcal C_0\setminus A,\ x\sim y\bigr\}.
\]
These quantities are $\mathcal R$-measurable. All edges between $K$ and $A\setminus K$ are closed, since $K$ is an open component of $A$. Let $J\coloneqq\{K\subset\mathcal C_B\}$. We claim that $J$ occurs if and only if at least one of the $m$ edges between $K$ and $\mathcal C_0\setminus A$ is open. If such an edge is open, then $K$ is connected to $B$, since $K$ is internally connected. Conversely, if $K$ is connected to $B$, consider an open path from $K$ to $B$ and its last vertex in $K$. The remaining part of the path avoids $K$ and starts with an edge from $K$ to a vertex $y$, which is not in $A$ and is connected to $B$ in the graph without $K$; thus $y\in\mathcal C_0\setminus A$. Since the $m$ edges are independent of $\mathcal R$, we obtain
\[
\mathbb P_p(J\,\lvert\,\mathcal R)=1-(1-p)^m.
\]
On $J^\c$, the set $K$ is not connected to $B$, and opening edges from $K$ to vertices outside $\mathcal C_0$ does not connect any further vertex to $B$. Hence $\mathcal C_B=\mathcal C_0$, $\mathcal N=N_0$ and $\mathcal M_i=m$ on $J^\c$. On $J$ we have $\mathcal C_B\supset\mathcal C_0\cup K$, and hence $\mathcal N\ge N_0+k$ and $\mathcal M_i\ge m$.

We split $J$ according to the size of $\mathcal M_i$. On $\{\mathcal M_i\ge k/(4p)\}$ we have $k\,4^{-\mathcal N}\le4p\,\mathcal M_i\,4^{-\mathcal N}$. On $J\cap\{\mathcal M_i<k/(4p)\}$ we have $m<k/(4p)$ and $4^{-\mathcal N}\le 4^{-(N_0+k)}$, so that
\[
\mathbb E_p\bigl[k\,4^{-\mathcal N}\1_J\1\{\mathcal M_i<k/(4p)\}\bigm\lvert\mathcal R\bigr]\le k\,4^{-(N_0+k)}\bigl(1-(1-p)^m\bigr)\1\{m<k/(4p)\}.
\]
The right-hand side vanishes if $m=0$. If $1\le m<k/(4p)$, we use $1-(1-p)^m\le pm$ and $-\log(1-p)\le2p$, which holds for $p\le1/2$. This gives $(1-p)^{-m}\le\e^{2pm}\le\e^{k/2}$ and
\[
k\,4^{-k}(1-p)^{-m}\le k\Bigl(\frac{\e^{1/2}}4\Bigr)^k\le k2^{-k}\le1.
\]
Consequently, the right-hand side of the previous display is at most
\[
p\,4^{-N_0}m(1-p)^m=p\,\mathbb E_p\bigl[4^{-\mathcal N}\mathcal M_i\1_{J^\c}\bigm\lvert\mathcal R\bigr].
\]
Combining both cases,
\[
\mathbb E_p\bigl[k\,4^{-\mathcal N}\1_J\bigm\lvert\mathcal R\bigr]\le4p\,\mathbb E_p\bigl[\mathcal M_i\,4^{-\mathcal N}\1_J\bigm\lvert\mathcal R\bigr]+p\,\mathbb E_p\bigl[\mathcal M_i\,4^{-\mathcal N}\1_{J^\c}\bigm\lvert\mathcal R\bigr]\le5p\,\mathbb E_p\bigl[\mathcal M_i\,4^{-\mathcal N}\bigm\lvert\mathcal R\bigr].
\]
Taking conditional expectations with respect to $\mathcal R_0$ yields~\eqref{eq:boundary-i}.
\end{proof}

\subsection{Reaching a set of vertices}\label{sec:target}

For a finite graph $H$ the contact process $(\xi^{H}_t)_{t\ge0}$ is a continuous-time Markov chain on the subsets of $V_H$ with generator
\[
\mathcal LF(A)=\lambda\sum_{x\in A}\sum_{y\in N_H(x)\setminus A}\bigl(F(A\cup\{y\})-F(A)\bigr)+\sum_{x\in A}\bigl(F(A\setminus\{x\})-F(A)\bigr),
\]
where $N_H(x)$ denotes the set of neighbours of $x$ in $H$.

\begin{proposition}\label{prop:target}
Let $H$ be a finite simple graph, $\emptyset\neq B\subset V_H$ and $o\in V_H\setminus B$. If $0<p\le1/2$ and $\lambda\ge20p$, then
\[
\mathbb P_\lambda\bigl(\xi^{H}_t\cap B\neq\emptyset\text{ for some }t\ge0\bigr)\ge\tfrac34\,\mathbb P_p(o\leftrightarrow B\text{ in }H).
\]
\end{proposition}

\begin{proof}
For $A\subset V_H$ let $F(A)\coloneqq1-\mathbb E_p[4^{-\#(A\cap\mathcal C_B)}]$. Then $F(\emptyset)=0$, $0\le F\le1$ and $F(\{o\})=(3/4)\,\mathbb P_p(o\in\mathcal C_B)=(3/4)\,\mathbb P_p(o\leftrightarrow B\text{ in }H)$. Let $A\subset V_H\setminus B$, and let $\mathcal N,\mathcal M$ be as in \cref{lem:boundary}. Since $\mathcal C_B$ does not depend on $A$, we have for $y\notin A$ and $x\in A$, respectively,
\[
\begin{aligned}
F(A\cup\{y\})-F(A)&=(3/4)\,\mathbb E_p\bigl[4^{-\mathcal N}\1\{y\in\mathcal C_B\}\bigr],\\
F(A\setminus\{x\})-F(A)&=-3\,\mathbb E_p\bigl[4^{-\mathcal N}\1\{x\in\mathcal C_B\}\bigr].
\end{aligned}
\]
Summing over the transitions of the contact process and using \cref{lem:boundary},
\[
\mathcal LF(A)=(3/4)\bigl(\lambda\,\mathbb E_p[\mathcal M\,4^{-\mathcal N}]-4\,\mathbb E_p[\mathcal N\,4^{-\mathcal N}]\bigr)\ge(3/4)(\lambda-20p)\,\mathbb E_p[\mathcal M\,4^{-\mathcal N}]\ge0.
\]
Let $\tau$ be the first time at which $\xi^{H}_t$ intersects $B$ or equals $\emptyset$. Since the state space is finite, $F(\xi^{H}_{t\wedge\tau})-\int_0^{t\wedge\tau}\mathcal LF(\xi^{H}_s)\dd s$ is a martingale, and hence $(F(\xi^{H}_{t\wedge\tau}))_{t\ge0}$ is a bounded submartingale. Moreover, $\tau<\infty$ almost surely. Indeed, in every time interval $[j,j+1]$, $j\in\mathbb N_0$, the event that there are no infection arrows in $H$ and every vertex of $H$ has a recovery mark has probability $\e^{-2\lambda\#E_H}(1-\e^{-1})^{\#V_H}>0$; these events are independent for different $j$, and each of them forces $\xi^{H}_{j+1}=\emptyset$. By bounded convergence,
\[
F(\{o\})\le\mathbb E_\lambda\bigl[F(\xi^{H}_\tau)\bigr]\le\mathbb P_\lambda\bigl(\xi^{H}_\tau\cap B\neq\emptyset\bigr),
\]
where we used $F(\emptyset)=0$ and $F\le1$.
\end{proof}

\subsection{From spatial spread to survival}\label{sec:exhaustion}

For $n\in\mathbb N$ let $H_n$ be the subgraph of $G$ induced by the ball $\{x\in V\colon\dist(o,x)\le n\}$ and let $S_n\coloneqq\{x\in V\colon\dist(o,x)=n\}$. Since $G$ is locally finite, $H_n$ is finite, and $S_n$ is the set of vertices of $H_n$ at maximal distance from $o$. We consider the contact process $\xi^{H_n}_t$ on $H_n$ started from $\{o\}$, constructed from the restriction of the graphical representation of $G$, and let
\[
\tau_n\coloneqq\inf\{t\ge0\colon\xi^{H_n}_t=\emptyset\},\qquad h_n\coloneqq\inf\{t\ge0\colon\xi^{H_n}_t\cap S_n\neq\emptyset\}.
\]
Thus $h_n$ is the first time at which the infection reaches distance $n$ from the root. Since $\xi^{H_n}_0=\{o\}$ and $o\notin S_n$, we have $h_n>0$, and $h_n<\tau_n$ means that the infection reaches distance $n$ before it dies out.

\begin{lemma}\label{lem:exhaustion}
Let $\lambda>0$ and $c>0$. If $\mathbb P_\lambda(h_n<\tau_n)\ge c$ for all $n\in\mathbb N$, then $\theta_{\mathsf{CP}}(\lambda)\ge c$.
\end{lemma}

\begin{proof}
Fix $u>0$ and integers $n\ge1$ and $K\ge2$, and let $\sigma_K$ be the first time at which $\xi^{H_n}_t$ contains at least $K$ vertices. On $\{\sigma_K\le\tau_n\le u\}$, each of the at least $K$ vertices that are infected at time $\sigma_K$ has a recovery mark in $(\sigma_K,u]$. By the strong Markov property at $\sigma_K$ and the independence of the recovery marks,
\[
\mathbb P_\lambda(\sigma_K\le\tau_n\le u)\le(1-\e^{-u})^K.
\]
On $\{h_n<\tau_n\le u,\ \sigma_K>\tau_n\}$, fewer than $K$ vertices are infected at any time before $\tau_n$, but at least $n+1$ distinct vertices have been infected at some time before $\tau_n$, since an infection path from $o$ to $S_n$ visits at least $n+1$ vertices, and each of them has recovered by time $\tau_n$. Thus at least $n$ recovery events of the finite Markov chain $\xi^{H_n}$ occur in $[0,u]$ before $\sigma_K$. Before $\sigma_K$, use a Poisson clock of rate $K$ to propose recoveries. If $m>0$ vertices are infected, accept each proposal with probability $m/K$ and, upon acceptance, recover one uniformly; reject proposals when $m=0$. The recovery count up to $u\wedge\sigma_K$ is then bounded by the proposal count in $[0,u]$, which has distribution $\mathrm{Poi}(Ku)$. Therefore,
\[
\mathbb P_\lambda(h_n<\tau_n\le u)\le(1-\e^{-u})^K+\mathbb P\bigl(\mathrm{Poi}(Ku)\ge n\bigr),
\]
and consequently, for all $n\ge1$ and $K\ge2$,
\[
\mathbb P_\lambda(\xi^{H_n}_u\neq\emptyset)=\mathbb P_\lambda(\tau_n>u)\ge c-(1-\e^{-u})^K-\mathbb P\bigl(\mathrm{Poi}(Ku)\ge n\bigr).
\]
By monotonicity in the graph, the events $\{\xi^{H_n}_u\neq\emptyset\}$ increase in $n$. Their union is $\{\xi_u\neq\emptyset\}$, since every infection path in $G$ uses finitely many arrows and is therefore contained in some $H_n$. Letting $n\to\infty$ for fixed $K$ yields $\mathbb P_\lambda(\xi_u\neq\emptyset)\ge c-(1-\e^{-u})^K$ for every $K\ge2$, and letting $K\to\infty$ gives $\mathbb P_\lambda(\xi_u\neq\emptyset)\ge c$ for every $u>0$. The events $\{\xi_u\neq\emptyset\}$ decrease in $u$, since every infection path to time $u$ has an initial segment reaching any earlier time, and their intersection over $u\in\mathbb N$ is the survival event. Hence $\theta_{\mathsf{CP}}(\lambda)\ge c$.
\end{proof}

\begin{proof}[Proof of \cref{thm:comparison}]
Let $n\in\mathbb N$. By \cref{prop:target} applied to $H=H_n$ and $B=S_n$,
\[
\mathbb P_\lambda(h_n<\tau_n)\ge\tfrac34\,\mathbb P_p(o\leftrightarrow S_n\text{ in }H_n)\ge\tfrac34\,\theta_{\mathsf{bond}}(p).
\]
The last inequality holds since an infinite open cluster of $o$ contains an open self-avoiding path from $o$ to $S_n$, whose initial segment up to its first visit to $S_n$ lies in $H_n$. If $\theta_{\mathsf{bond}}(p)=0$, there is nothing to prove; otherwise \cref{lem:exhaustion} with $c=(3/4)\,\theta_{\mathsf{bond}}(p)$ completes the proof.
\end{proof}

\section{Proof of \cref*{prop:tree}}\label{sec:proof-tree}

Throughout this section, $T$ is an infinite, locally finite tree with root $o$. For $x\neq o$ we denote the parent of $x$ by $x^-$, and we write $x\le y$ if $x$ lies on the path from $o$ to $y$. For vertices $x,y$ we let $x\wedge y$ be the vertex of largest distance to $o$ with $x\wedge y\le x$ and $x\wedge y\le y$. We use the graphical representation of \cref{sec:setting} on $T$.

\subsection{Quasi-independent percolation}\label{sec:sitebond}

The following second-moment bound is a special case of Lyons' result on quasi-independent percolation~\cite{Lyons1989}; we include the short proof for completeness.

\begin{lemma}\label{lem:lyons}
Let $\cK\subset V$ be a random set with $o\in\cK$ such that $y\in\cK$ implies $y^-\in\cK$ for all $y\neq o$. Assume that there are $q\in(0,1)$ and $M\ge1$ such that, for all vertices $x,y$,
\begin{equation}\label{eq:quasi}
\mathbb P(x\in\cK)=q^{|x|}\qquad\text{and}\qquad\mathbb P(x\in\cK,\ y\in\cK)\le M\,q^{|x|+|y|-|x\wedge y|}.
\end{equation}
If $q\br(T)>1$, then $\mathbb P(\#\cK=\infty)>0$.
\end{lemma}

\begin{proof}
Since $1/q<\br(T)$, the random walk on $T$ with conductances $q^{|e|}$ is transient~\cite{Lyons1990}, see also~\cite[Chapter~3]{LyonsPeres2016}. Hence there is a unit flow $\theta$ from $o$ to infinity with finite energy $\mathcal E\coloneqq\sum_e\theta(e)^2q^{-|e|}<\infty$~\cite[Chapter~2]{LyonsPeres2016}. We take for $\theta$ the unit current flow, which is non-negative on trees, since the normalised voltage $x\mapsto\mathbb P_x(\text{the walk ever hits }o)$ is non-increasing along every path away from $o$.
For a vertex $x\neq o$ we write $\theta(x)$ for the flow through the edge $\{x^-,x\}$. By flow conservation, $\sum_{|x|=n}\theta(x)=1$ for all $n\in\mathbb N$, and $\sum_{|x|=n,\,x\ge y}\theta(x)=\theta(y)$ for all $y\neq o$ with $|y|\le n$.

Fix $n\in\mathbb N$ and set $X_n\coloneqq\sum_{|x|=n}\theta(x)q^{-n}\1\{x\in\cK\}$. Then $\mathbb E[X_n]=1$, and by~\eqref{eq:quasi},
\[
\mathbb E[X_n^2]\le M\sum_{|x|=|y|=n}\theta(x)\theta(y)q^{-|x\wedge y|}.
\]
For every vertex $w$ we write $q^{-|w|}=1+\sum_{e}(q^{-|e|}-q^{-|e|+1})$, where the sum runs over the edges $e$ on the path from $o$ to $w$. An edge $e=\{y^-,y\}$ lies on the path from $o$ to $x\wedge y'$ if and only if $x\ge y$ and $y'\ge y$. Therefore,
\[
\sum_{|x|=|y'|=n}\theta(x)\theta(y')q^{-|x\wedge y'|}=1+\sum_{|e|\le n}\bigl(q^{-|e|}-q^{-|e|+1}\bigr)\theta(e)^2\le1+\mathcal E.
\]
By the Paley--Zygmund inequality,
\[
\mathbb P(\cK\cap\{x\colon|x|=n\}\neq\emptyset)\ge\mathbb P(X_n>0)\ge\frac{\mathbb E[X_n]^2}{\mathbb E[X_n^2]}\ge\frac1{M(1+\mathcal E)}.
\]
Since $\cK$ is a subtree containing $o$ and $T$ is locally finite, $\cK$ is infinite if and only if it intersects every level. The events on the left-hand side decrease in $n$, and hence $\mathbb P(\#\cK=\infty)\ge1/(M(1+\mathcal E))>0$.
\end{proof}

\subsection{The single-infection process on trees}\label{sec:sir-tree}

We construct the single-infection process of \cref{sec:chain} on $T$ from the graphical representation. Let $\sigma_o\coloneqq0$. If $\sigma_x$ is defined and $y$ is a child of $x$, let $\rho_x\coloneqq\inf(D_x\cap(\sigma_x,\infty))$ be the first recovery mark at $x$ after $\sigma_x$ and set
\[
\sigma_y\coloneqq\begin{cases}\inf\bigl(A_{(x,y)}\cap(\sigma_x,\rho_x)\bigr)&\text{if }\sigma_x<\infty,\\ \infty&\text{otherwise,}\end{cases}
\]
with $\inf\emptyset\coloneqq\infty$. Thus $x$ is infected at time $\sigma_x$, stays infected until $\rho_x$, and infects a child $y$ at the first arrow from $x$ to $y$ during this period. Let $\cK\coloneqq\{x\in V\colon\sigma_x<\infty\}$.

\begin{lemma}\label{lem:sir-tree}
Let $\lambda>0$, $q\coloneqq\lambda/(1+\lambda)$ and $M\coloneqq2(1+\lambda)/(1+2\lambda)\le2$.
\begin{enumerate}[label=(\alph*)]
\item\label{it:sir-law} The set $\cK$ has the same distribution as the set $\cS_\lambda$ of \cref{sec:chain} on $T$, and it satisfies~\eqref{eq:quasi} with these $q$ and $M$.
\item\label{it:sir-cp} Let $x\in V$ with $\sigma_x<\infty$, and let $H$ be a subgraph of $T$ containing the path from $o$ to $x$, for instance a ball $H_n$ with $n\ge|x|$. Then $x$ is infected at time $\sigma_x$ in the contact process on $H$ started from $\{o\}$, that is, $x\in\xi^H_{\sigma_x}$.
\end{enumerate}
\end{lemma}

\begin{proof}
For part~\ref{it:sir-law}, let $T_x$ be the time from $\sigma_x$ to the next recovery mark at $x$, and $E_y$ the time from $\sigma_{y^-}$ to the next arrow from $y^-$ to $y$ (with each starting time replaced by $0$ when it is infinite). Then $y\in\cK$ if and only if $y^-\in\cK$ and $E_y<T_{y^-}$. The time $\sigma_x$ only depends on the processes attached to strict ancestors of $x$, which are independent of $D_x$ and the outgoing arrow processes $A_{(x,y)}$. Conditional on these ancestral processes, the futures of $D_x$ and the outgoing arrow processes after $\sigma_x$ (or after $0$ if $\sigma_x=\infty$), shifted back to start at $0$, are independent Poisson processes with their original rates. Proceeding generation by generation, we obtain that all $T_x$ and $E_y$ are independent and exponentially distributed with parameters $1$ and $\lambda$, respectively. On a tree, a vertex is reached from $o$ by a directed path of open edges if and only if all edges from parent to child on the path from $o$ to this vertex are open; hence $\cK$ has the distribution of $\cS_\lambda$.

Every edge on the path from $o$ to $x$ starts at a different vertex, and $\mathbb P(E_y<T_{y^-})=\lambda/(1+\lambda)=q$. Hence $\mathbb P(x\in\cK)=q^{|x|}$. Let $w\coloneqq x\wedge y$. If $w\in\{x,y\}$, say $w=x$, then $\mathbb P(x,y\in\cK)=q^{|y|}=q^{|x|+|y|-|w|}$. Otherwise the paths from $w$ to $x$ and to $y$ leave $w$ through different children $a\neq a'$ of $w$. The edges $\{w,a\}$ and $\{w,a'\}$ are both open with probability
\[
\mathbb E\bigl[(1-\e^{-\lambda T_w})^2\bigr]=1-\frac2{1+\lambda}+\frac1{1+2\lambda}=\frac{2\lambda^2}{(1+\lambda)(1+2\lambda)}=M\,q^2,
\]
and all other edges on the two paths start at different vertices. Therefore $\mathbb P(x,y\in\cK)=M\,q^{|x|+|y|-|w|}$, which proves~\eqref{eq:quasi}.

Part~\ref{it:sir-cp} holds by construction: the arrows at the times $\sigma_y$, $o<y\le x$, together with the absence of recovery marks at $y^-$ during $(\sigma_{y^-},\sigma_y]$, form an infection path from $(o,0)$ to $(x,\sigma_x)$ inside $H$.
\end{proof}

\begin{proof}[Proof of \cref{prop:tree}]
Let $\lambda>\pc(T)/(1-\pc(T))$. Then $q\coloneqq\lambda/(1+\lambda)>\pc(T)=1/\br(T)$ by~\eqref{eq:lyons}, that is, $q\br(T)>1$. By \cref{lem:lyons,lem:sir-tree}, the probability $c\coloneqq\mathbb P(\#\cK=\infty)$ is positive. Let $n\in\mathbb N$ and let $H_n$, $S_n$, $h_n$ and $\tau_n$ be as in \cref{sec:exhaustion}, with $G=T$. If $\cK$ is infinite, it contains a vertex $x$ with $|x|=n$, and by \cref{lem:sir-tree}\ref{it:sir-cp} we have $x\in\xi^{H_n}_{\sigma_x}\cap S_n$. Hence $h_n\le\sigma_x<\infty$ and $h_n<\tau_n$, so that $\mathbb P_\lambda(h_n<\tau_n)\ge c$ for all $n\in\mathbb N$. By \cref{lem:exhaustion}, $\theta_{\mathsf{CP}}(\lambda)\ge c>0$, and hence $\lc(T)\le\lambda$.
\end{proof}

\section{Proof of \cref*{prop:chain}}\label{sec:proof-chain}

\begin{proof}[Proof of \cref{prop:chain}\ref{it:chain-perc}]
Let $p\in(0,1]$. The edges of the paths between $v_k$ and $v_{k+1}$ are disjoint for different $k$, and $v_k$ is connected to $v_{k+1}$ by an open path inside the $k$-th block with probability $1-(1-p^2)^k$. Hence
\[
\mathbb P_p(o\leftrightarrow\infty)\ge\prod_{k\in\mathbb N}\bigl(1-(1-p^2)^k\bigr)>0,
\]
since all factors are positive and $\sum_k(1-p^2)^k<\infty$. Thus $\pc(G_{\mathsf{chain}})=0$.

Let $T$ be an infinite subtree of $G_{\mathsf{chain}}$. Every vertex $v_k$ separates the finite set of vertices of the blocks $1,\dots,k-1$ from the rest of the graph. Hence every infinite self-avoiding path in $T$ visits $v_k$ for all large $k$. If $T$ contained two disjoint rays, they would both pass through some $v_k$, which is impossible. Thus $T$ consists of a single ray together with finite subtrees, and every cutset consisting of a single edge of this ray far away from the root shows that $\br(T)\le1$. By~\eqref{eq:lyons}, $\pc(T)=1$.
\end{proof}

\begin{proof}[Proof of \cref{prop:chain}\ref{it:chain-sir}]
For $k\in\mathbb N$ let $F_k$ be the event that there is $1\le i\le k$ such that both directed edges $(v_k,w_{k,i})$ and $(w_{k,i},v_{k+1})$ are open. The only neighbours of $v_{k+1}$ in blocks $1,\dots,k$ are the vertices $w_{k,i}$, whose only other neighbour is $v_k$. Hence a directed open path from $o$ that reaches $v_{k+1}$ for the first time uses one of the pairs of edges in $F_k$, and $\{v_{k+1}\in\cS_\lambda\}\subset\bigcap_{j\le k}F_j$. If $\cS_\lambda$ is infinite, then $v_{k}\in\cS_\lambda$ for all $k$, since otherwise $\cS_\lambda$ is contained in the finite set of vertices of the blocks before $v_k$. The events $F_k$ depend on disjoint sets of the variables $T_x$ and $E_{(x,y)}$ and are therefore independent.

Let $q\coloneqq\lambda/(1+\lambda)=\mathbb P(E_{(w,v)}<T_w)$. Conditionally on $T_{v_k}$, the $k$ events in the definition of $F_k$ are independent with probability $(1-\e^{-\lambda T_{v_k}})q$ each. Hence
\[
\mathbb P(F_k^\c)=\mathbb E\Bigl[\bigl(1-q(1-\e^{-\lambda T_{v_k}})\bigr)^k\Bigr]\ge\mathbb P\bigl(T_{v_k}\le1/(\lambda k)\bigr)\Bigl(1-\frac1k\Bigr)^k\ge\frac1{8\lambda k}
\]
for all $k\ge\max\{2,1/\lambda\}$, where we used $1-\e^{-\lambda T}\le\lambda T$, $q\le1$, $(1-1/k)^k\ge1/4$ for $k\ge2$ and $1-\e^{-u}\ge u/2$ for $u\in[0,1]$. Since $\sum_k1/k=\infty$, we obtain $\mathbb P(\#\cS_\lambda=\infty)\le\prod_{k}\mathbb P(F_k)=0$.
\end{proof}

\begin{proof}[Proof of \cref{prop:chain}\ref{it:chain-cp}]
This follows from part~\ref{it:chain-perc} and \cref{cor:comparison}. Alternatively, for every $K\in\mathbb N$, delete the edges $\{w_{k,i},v_{k+1}\}$ for $k>K$ and $2\le i\le K+1$. In the resulting subgraph, the vertices $v_k$, $k>K$, as centres, the vertices $w_{k,1}$ as connecting paths and the leaves $w_{k,2},\dots,w_{k,K+1}$ attached to $v_k$ form a copy of $H_{K,1}$. Graph monotonicity and \cref{sec:converse} then apply.
\end{proof}

\section*{AI disclosure}
This note was produced with substantial assistance from the AI systems Claude (Anthropic) and ChatGPT/Codex (OpenAI). The authors have checked all statements and proofs and take full responsibility for the content of this note.

\section*{Acknowledgement}
We thank Daniel Valesin for his question on an earlier version of~\cite{JahnelLuechtrathMoench2025}, which motivated this note. BJ and LL acknowledge the financial support of the Leibniz Association within the Leibniz Junior Research Group on \emph{Probabilistic Methods for Dynamic Communication Networks} as part of the Leibniz Competition as well as the Berlin Cluster of Excellence MATH+ through the project EF45-3 on \emph{Data Transmission in Dynamical Random Networks}.

{\small
\bibliographystyle{amsalpha}
\bibliography{references}
}

\end{document}